\documentclass[12pt]{amsart}
\usepackage{latexsym,amssymb,amsmath}
\usepackage[all,ps,cmtip]{xy} 

\newtheorem{theo}{Theorem}
\newtheorem{coro}{Corollary}
\newtheorem{prop}{Proposition}
\newtheorem{lemm}{Lemma}

\def\fg{\mathfrak{g}}
\def\lra{\longrightarrow}
\def\cO{\mathcal{O}}
\def\PP{\mathbf{P}}\def\AA{\mathbf{A}}
\def\RR{\mathbf{R}}\def\CC{\mathbf{C}}\def\HH{\mathbf{H}}\def\OO{\mathbf{O}}

\begin{document}

\title{Hyperk\"ahler manifolds from \\ the Tits-Freudenthal magic square}
\author{Atanas Iliev, Laurent Manivel}


\begin{abstract}
We suggest a way to associate to each  Lie algebra of type $G_2$, $D_4$, $F_4$, $E_6$, $E_7$, $E_8$ a family of 
polarized hyperk\"ahler fourfolds, constructed as parametrizing certain families of cycles of hyperplane  
sections of certain homogeneous or quasi-homogeneous varieties. These cycles are modeled on the 
Legendrian varieties studied by Freudenthal in his geometric approach to the celebrated Tits-Freudenthal 
magic square of Lie algebras.  
\end{abstract}

\maketitle

\section{Introduction}

Hyperk\"ahler manifolds, also called irreducible compact holomorphic symplectic varieties, are a notoriously
mysterious and challenging class of varieties. A very frustrating question is to decide whether there could 
exist many examples besides deformations of Hilbert schemes $Hilb_n(S)$ of points on a K3 surface $S$, 
and the closely related generalized Kummer varieties. Even the general deformation of $Hilb_n(S)$
is not known. If we had a polarization and consider deformations of polarized holomorphic symplectic varieties, 
the general deformation has only been described in a restricted number of cases, each exhibiting some 
very special and interesting geometry. 

In this paper, we stress an unexpected relation between the exceptional groups and certain families 
of polarized holomorphic symplectic varieties. Briefly, this goes as follows. To each simple complex Lie
algebra one can associate a particular projective homogeneous variety $X_{ad}$, which we call the adjoint variety, 
whose automorphism group is the adjoint form of the Lie group. The adjoint variety is covered by a 
family of special subvarieties $C$ that we call Legendrian cycles. Then our key observation is that for each exceptional
group, there exists a projective variety $X$ of the same dimension as $X_{ad}$, and covered by the same 
type of special subvarieties, i.e. by Legendrian cycles. Moreover there is a suitable cycle space $P$ for these Legendrian cycles plus certain of their 
degenerations. And the crucial observation is that for a generic hyperplane section $X_H$ 
of $X$ the subspace $P_H$ of $P$ parametrizing the Legendrian cycles $C$ which are contained in $X_H$ 
is a  hyperk\"ahler fourfold. 

This has to be taken with a grain of salt. Indeed when we start from $F_4$, the resulting variety $P_H$ 
is only a blow-up of a   hyperk\"ahler fourfold at nine points. Moreover these fourfolds are essentially the
same as those that are obtained from $E_6$. Starting from $E_7$, the related geometry becomes 
much more complicate, and in the cases $E_7$ and $E_8$ we had not been able to achieve the 
expected construction that should follow the above line of the argument. 
Nevertheless, we believe that its existence also in these two cases remains an interesting challenge. 

\medskip
The following table summarizes the  hyperk\"ahler geometries we have encountered. \medskip

\begin{center}
\begin{tabular}{|c|l|}\hline  
$G_2$ & Fano variety of lines on a cubic fourfold \cite{bd} \\
$D_4$ & Hilbert square of determinantal quartic surfaces \\
$F_4, E_6$ & Zero-locus of a 3-form in $10$ variables \cite{dv} \\
$E_7, E_8$ & ?? \\ \hline
\end{tabular}
\end{center}

\medskip\noindent {\it Structutre of the paper}.

In Section 2 are collected some facts about the exceptional series of Lie 
groups, the Tits-Freudenthal magic square and the associated geometries, 
which are related to our study.  

Section 3 provides details on the construction we have just sketched. 

Section 4 treats the case of $D_4$, for which triality 
plays an interesting role. 

Section 5 discusses the case of $F_4$, for which we construct a suitable 
cycle space for the Legendrian cycle spaces and their degenerations. 

In Section 6 we construct the associated  hyperk\"ahler fourfold and we show, 
following a suggestion of Kuznetsov, that it is isomorphic to the fourfold obtained form $E_6$. 


\medskip\noindent {\it Acknowledgements.} We thank A. Kuznetsov, D. Orlov, K. O'Grady, G. Kaputska
and M. Kaputska for their useful hints and comments. 

This work has been carried out in the framework of the Labex Archimede (ANR-11-LABX-0033) 
and of the A*MIDEX project (ANR-11-IDEX-0001-02), funded by the ``Investissements d'Avenir" French 
Government programme managed by the French National Research Agency, 
and the Seoul National University grant 0450-20130016.

\section{The exceptional series and the magic square}\label{section 2}

Let $G$ be a simple complex Lie group, and let $\fg$ denote its Lie algebra. The action of $G$ on $\PP\fg$
has a unique closed orbit, which we call the {\it  adjoint variety}. It is the projectivization of the minimal non trivial nilpotent 
orbit in $\fg$, and the Kostant-Kirillov-Souriau symplectic structure on this nilpotent orbit descends to a contact 
structure on the adjoint variety. 

Choose a point $x$ of the adjoint variety $X$, and denote by $P$ the stabilizer of $x$ in $G$. Then $P$ is a parabolic subgroup
of $G$, and there is an induced isomorphism $X\simeq G/P$. The isotropy representation of $P$ on $T_xX$ has a unique invariant 
hyperplane $V$ (the contact hyperplane), which is (except in type A) an irreducible representation of a given 
Levi factor of the parabolic group $P$. Recall that a Levi factor is a reductive subgroup of $P$, such that the projection
to $P/P^u$, the quotient of $P$ by its unipotent radical, is an isomorphism; it is uniquely defined up to conjugation.

 Let us denote by $H$ the semisimple part of the chosen Levi factor. Geometrically, the set $C$  
of projective lines in $\PP\fg$ contained in $X$ and passing through $x$  is acted on transitively by $H$; its linear span in 
$T_xX$ is the contact hyperplane $V$.  

The reduction of $X\subset \PP\fg$ to $C\subset \PP V$ is a general process
which has been studied in detail in \cite{lm}--\cite{lm4}. 
One of the important aspects of this study is the observation that this reduction process provides a geometric
transition between the lines of the Tits-Freudenthal magic square, recently enlarged to a rectangle by Deligne and
Gross. Starting from the adjoint varieties of the simple Lie algebras in the Deligne's extended exceptional series, one gets 
the following list of groups and varieties: 

$$\begin{array}{cccccccc} 
G &  & G_2 & SO_8 & F_4  & E_6  & E_7  & E_8 \\
h_G^\vee &&4&6&9 & 12 & 18 & 30 \\
 a &&-\frac{2}{3}&0&1&2&4&8 \\
H &  & SL_2 &   SL_2^3 & Sp_6 & SL_6 & Spin_{12} & E_7 \\
C & & v_3\PP^1 & (\PP^1)^3& IG(3,6) & G(3,6) & S_{12} & Fr 
\end{array}$$

\medskip
We include in this table the dual Coxeter number $h_G^\vee$, and the closely related parameter $a=h_G^\vee/3-2$.
The dual Coxeter number is essentially 
the parameter used by Deligne in \cite{del}. The parameter $a$ makes more visible the connection with 
the Tits-Freudenthal magic square, where the Lie algebras of type $F_4, E_6, E_7, E_8$ are constructed from the 
pairs $(\AA, \OO)$, with $\AA=\RR,\CC,\HH,\OO$ a normed algebra of real dimension $a=1,2,4,8$ -- see \cite{lm} and 
the references therein. In particular
the existence of the exceptional series appears as a byproduct of the existence of the Cayley (or rather Cayley-Salmon)
algebra of octonions. In our study we will meet again the name of Cayley, in an seemingly 
completely independent, but in fact implicitly connected context: the subgroups 
$H$ from the series $Sp_6,  SL_6, Spin_{12},E_7$ reappear on the third line of the magic square, 
where they can be reconstructed from the pairs $(\AA, \HH)$. 

A variant of the original Tits-Freudenthal construction, called the triality construction, 
shows that it is natural to include $SO_8$ in the exceptional
series with parameter $a=0$, see \cite{lm} and \cite{del,delgross}. 
From this perspective, the group   
$Spin_8$ and the triality are at the origin of the exceptional series. 
Finally, $G_2$ is deduced from $SO_8$ (or rather by $Spin_8$) by a triple folding,
and the denominator of the corresponding parameter $a=-2/3$ can be thought 
of as a reminiscence of this folding.
The dimensions of $X$ and $V$ as above, and of $C$ from the above table, are connected by the 
following dimension formulas:
$$\dim C=3a+3, \qquad \dim V = 6a+8=\dim X-1.$$
The fact that the cone $\hat{C}$  over $C$ has half the dimension of $V$ is not a coincidence. 
There exists an $H$-invariant symplectic structure on $V$ with respect to which 
$\hat{C}$ is Lagrangian \cite{lm1}; and $C$, which by itself is a Fano manifold of index $2a+2$, 
is said to be Legendrian, see \cite{lm4}, \cite{buc} for a systematic study of Legendrian varieties.


\section{Legendrian cycles and symplectic manifolds}\label{legendrian}

We suggest a procedure that allows to construct families of holomorphic symplectic manifolds, starting from an
adjoint variety $X_{ad}$ of the exceptional series. 

\medskip\noindent {\bf First step}. Replace $X_{ad}$ by a rational variety $X$ of the same dimension $6a+9$,
covered by cycles of type $C$. This is actually a peculiar step of the construction. 

The variety $X$ is again rational and Fano of index $3a+4$, different from the index of $X_{ad}$, 
which by itself is $3a+5$; cf. the tables in \S 2 of \cite{cp} or in \cite{snow} 
for the invariants for the adjoint varieties of the classical Lie groups.  
We will denote by $L$ the (very ample) line bundle on $X$ such that $-K_X=(3a+4)L.$

\medskip\noindent {\bf Second step}.
 Find a good parameter space $P$ for cycles of type $C$ in $X$ and their degenerations,
which we call Legendrian cycles in $X$. In particular $P$ must be smooth and endowed 
with a family of cycles ${\mathcal C}\subset P\times X$.

$$\begin{array}{cccccccc} 
G &  & G_2 & SO_8 & F_4  & E_6  & E_7  & E_8 \\
X && \PP^5 &  (\PP^3)^3& IG(3,9) & G(3,10) & OG(6,15) & ?? \\
P && G(2,6) & G(2,4)^3 &Bl G(6,9) & G(6,10) & Bl G(12,15) & ?? 
\end{array}$$

\medskip
Note that although $P$ is close to be homogeneous, it is not always acted on transitively by the automorphism group of 
$X$. Two non transitive actions appear for $G=F_4$ and $E_7$, where $Bl$ means that some blowing-up is required. 
Moreover in the first of these two cases, 
the automorphism group of $X_{Leg}$ is the odd symplectic group $Sp(9)$, which is not semisimple or even reductive. 
This will be a source of additional technical complications.

The dimension of $P$ is given by the simple formula
$$\dim P=6a+12.$$ 

\noindent {\bf Third step}. Let $X_H=X\cap H$ be a general hyperplane section of $X$, defined by a general element
$h\in H^0(X,L)$. In particular $X_H$ is 
smooth of dimension $6a+8$. We will check that the Hodge cohomology group
$H^{3a+5,3a+3}(X_H)$
has dimension \underline{one}, i.e.  
 $$H^{3a+5,3a+3}(X_H)=\CC \Omega$$
for some closed non-exact $(3a+5,3a+3)$-form $\Omega$ on $X_H$. 
Moreover $H^{p,q}(X_H)=0$ for $p+q=6a+8$ and $q<3a+3$. In other words, the middle
dimensional Hodge structure of $X_H$ is {\it of $K3$ type}. 

\medskip\noindent {\bf Fourth step}. Consider the space $P_H\subset P$ of Legendrian cycles $C$ that are contained in $H$. 
The parameter space $P$ will be constructed in such a way that $E=p_*q^*L$ will be a rank $6a+8$ 
vector bundle on $P$ with $H^0(P,E)\cong H^0(X,L).$
In particular, $h$ will define a section $h_E$ of $E$ with zero-locus $P_H=Z(h_E)$,  
 which turns out to be smooth and of dimension four. 

 \medskip\noindent {\bf Fifth step}. Consider the restriction of the family of Legendrian cycles ${\mathcal C}\subset 
P_H\times C_H$, with its two projections $p_H$ and $q_H$. The Abel-Jacobi map 
$$p_{H*}q_H^* : H^{3a+5,3a+3}(X_H)\longrightarrow H^{2,0}(P_H)$$
maps $\Omega$ to a holomorphic two-form $\omega$ on $P_H$  which is generically non-degenerate. 

 \medskip\noindent {\bf Sixth step}. Find, if necessary, a suitable modification $Z_H$ of $P_H$ with trivial 
canonical bundle, and conclude that $Z_H$ is  hyperk\"ahler. 

\bigskip\noindent {\it Comments}. For $a=-2/3, 0, 2$ the sixth and last step is not necessary: we get directly
a  hyperk\"ahler structure on $P_H$, and no birational modification is needed. Moreover in these cases $P$ appears
to be homogeneous. In particular the cycles in $X$ isomorphic to $C$ form a complete family, 
and there is no need to take degenerations into account.

The hyperk\"ahler vaieties from the case $a=2$ 
have already been discovered by Debarre and Voisin in \cite{dv}, 
where they prove that the result of virtually the same construction is 
a locally complete family of polarized hyperk\"ahler fourfolds, 
which are deformations of Hilbert squares of K3 surfaces. 

 In fact, the case $a=-2/3$ is even more classical:
it is nothing else than the  hyperk\"ahler structure on the Fano variety of lines of a cubic fourfold discovered by Beauville 
and Donagi in \cite{bd}. 

From our perspective, the special structures of these two very classical varieties have their hidden 
origin respectively in the exceptional groups $E_6$ and $G_2$. 

From the triality perspective, the case $a=0$ is the most important. It could have been discovered much earlier,
being defined by a very simple and natural vector bundle on a product of three four-dimensional quadrics. 
It turns out that this construction does not yield a  locally complete family of projective  hyperk\"ahler manifolds.
In fact the symplectic varieties that we obtain are Hilbert squares of determinantal quartic surfaces in 
$\PP^3$. 

The case $a=1$ will be the object of section \ref{a=1}. As already mentioned, a new feature here is the necessity 
to take into account degenerations of our Legendrian cycles $C$. The cycle space $P$ will be defined as a blow-up 
of a Grassmannian. As explained in the fourth step above, it is endowed with a vector bundle, but the  zero-locus
of its general section is not a  hyperk\"ahler manifold. In fact it  contains $9$ special copies of $\PP^3$, and it is
only after contracting these that we finally obtain a  hyperk\"ahler manifold. Alternatively, one can consider, 
directly on the Grassmannian, the zero-locus of a sheaf which is not locally free. 
%

In the case $a=4$  there is an additional complication: we would expect $X$ to be Fano of index  $10$, 
but it is only of index $5$. The reason for this 
discrepancy is that our Legendrian cycles $C$ in $X$ are copies of $OG_+(6,12)$ inside $OG(6,15)$;
and the Pl\"ucker polarization on $OG(6,15)$, which is certainly primitive, becomes divisible by two when 
restricted to $OG_+(6,12)$, a phenomenon which is directly related to the existence of the Spin 
representations. So one should modify the construction and ask that the cycle family $\mathcal{C}$ 
be endowed with a line bundle $\mathcal{L}$ such that $\mathcal{L}^2=p^*O(1)$, and $E=q_*L$ be 
a vector bundle of rank $32$ (given at the generic point of $P$ by a half-spin representation of $Spin_{12}$).
This is something we have not managed yet.

The case $a=8$ is the most mysterious, and in this final case we have 
no suggestion for eventual constructions of $X$ and $P$,
which should follow the conjectured Steps 1-6, as stated above. 
 

An intermediate case, first noticed by Deligne as a complement to the exceptional series
is that of $a=6$. Algebraically, one can complete
the magic square with the help of  the intermediate algebra of sextonions, 
which is normed for a degenerate quadratic form. 
Geometrically, we know that the geometries associated with the Freudenthal square can still be 
constructed from the sextonions, but they yield singular varieties. 
Therefore it is plausible that our construction may be 
used for a search of new singular hyperk\"ahler varieties of dimension 4. 
If this is indeed the case, the next question will obviously be that of 
the existence of symplectic resolutions of these fourfolds. 


\section{Determinantal quartics and triality}

In this and the next sections we keep the notation from Section \ref{legendrian},
and start with the case $a=0$. The Legendrian cycle corresponding to this case is 
$C=\PP^1\times \PP^1\times \PP^1\subset \PP^7$, see e.g. Table I on p. 89 of \cite{buc}. 
Its dual or tangent variety is the quartic 
hypersurface whose equation is Cayley's famous $2\times 2\times 2$ hyperdeterminant \cite{gkz}. 

The adjoint variety of $D_4$ is the orthogonal Grassmannian $OG(2,8)$, 
which has dimension $9 = 6a+9$ and index $5 = 3a+5$, see \S 2 in \cite{cp}.  
As already discussed in Section \ref{legendrian}, 
we replace the 9-dimensional adjoint variety of $D_4$, which is 
the orthogonal grassmannian $X_{ad}=OG(2,8)$, by  $X=\PP^3\times \PP^3\times \PP^3$. This variety is covered by copies of
$C$ parametrized by $P=G(2,4)\times G(2,4)\times G(2,4)$, a triple product of 4-dimensional
quadrics. To be more specific we will introduce  four-dimensional vector spaces $V_1, V_2, V_3$ and
let $X=\PP V_1\times \PP V_2\times \PP V_3$. 
$X$ is evidently a Fano manifold of dimension $9 = 6a+9$ and index 
$4 = 3a+4$: the anticanonical divisor $-K_X=4L$ for $L=\mathcal{O}_X(1,1,1)$. 

The parameter space for the Legendrian cycles $C \subset X$ is 
$P=G(2,V_1)\times G(2,V_2)\times G(2,V_3)$, and $P$ has dimension $6a+12 = 12$. 

We denote by  $T_i$ and $Q_i$ the rank two tautological bundles and quotient bundles  
on the Grassmannians $G(2,V_i)$, $i = 1,2,3$. The vector bundle $E$ on $P$ of 
rank $6a+8 = 8$ is $E=T_1^*\otimes T_2^*\otimes T_3^*$. 

Consider a general element of  
$$H^0(X,L)=H^0(P,E)= V_1^*\otimes V_2^*\otimes V_3^*,$$
considered either as a section $h$ of $L$, whose zero-locus is a hyperplane section $X_H$ of $X$, 
or as a section $h_E$ of $E$, whose zero locus in $P$ is a four-dimensional subvariety $P_H$. 
Note that  since $\det (E)=\cO(4,4,4)$, the canonical bundle of $P_H$ is trivial.

\begin{theo}
 $P_H$ is an irreducible holomorphic symplectic manifold.
\end{theo}

One way to show this is to break the symmetry of $h\in  V_1^*\otimes V_2^*\otimes V_3^*$, and consider
it as a morphism $h_3 : V_3\lra  V_1^*\otimes V_2^*$ -- or equivalently, a four dimensional space of 
sections of $\cO(1,1)$ on $\PP(V_1)\times \PP(V_2)$. Generically, such a linear system defines a K3
surface $S_3\subset \PP(V_1)\times \PP(V_2)$, parametrizing pairs of lines $(\ell_1,\ell_2)$ such that 
 $h$ vanishes on $\ell_1\otimes\ell_2\otimes V_3$. 

In particular $\ell_1$ belongs to the projection of $S_3$ to  $\PP(V_1)$ if and only if the image of 
$\ell_1$ by   $h_1: V_1\lra  V_2^*\otimes V_3^*\simeq Hom(V_2,V_3^*)$ is generated by a non injective 
morphism. For $h$ general the rank drops by two only in codimension four, hence nowhere 
since we are on $\PP^3$. There the projection of $S_3$ to  $\PP(V_1)$ is in fact an 
isomorphism to a determinantal quartic surface. 

Note also that  $Hom(V_2,V_3^*)\simeq  Hom(V_3,V_2^*)$, and that this isomorphism preserves the 
rank. This implies that if $h$ vanishes on $\ell_1\otimes\ell_2\otimes V_3$,  then there exists a 
line $\ell_3$ such that  $h$ vanishes on $\ell_1\otimes V_2\otimes \ell_3$. Moreover, this line 
is unique if $h$ is general. As a consequence, all the three surfaces $S_1, S_2, S_3$ are  
isomorphic to the same K3 surface $S$, of which we get three in general different representations 
as a determinantal quartic surface: 
$$
Q_1\subset \PP(V_1),  \ \ \ Q_2\subset \PP(V_2),  \ \ \ Q_3\subset \PP(V_3)
.$$

It is remarkable that this has already been known to Arthur Cayley  (see \cite{cayley} 
and \cite{fggl} for more historical information about this classical construction). 
Cayley has also observed that the  correspondence 
$\ell_1\otimes\ell_2\otimes V_3\mapsto\ell_1\otimes V_2\otimes \ell_3$ can be iterated, and that the iteration 
$$\ell_1\otimes\ell_2\otimes V_3\mapsto\ell_1\otimes V_2\otimes \ell_3\mapsto V_1\otimes \ell '_2\otimes \ell_3
\mapsto\ell '_1\otimes \ell '_2\otimes V_3$$
defines a non-trivial automorphism of $S$. 
By its construction this automorphism is a byproduct of 
the triality for $SO(8)$, and so one can call it the triality automorphism. 

\medskip Theorem 1 will be a consequence of the following more precise statement. 
 
\begin{prop}
$P_H$ is isomorphic to $Hilb^2(S)$.
\end{prop}
 
\proof A point of $P_H$ is a triple $(T_1,T_2,T_3)$ in $G(2,V_1)\times G(2,V_2)\times G(2,V_3)$ such
that $h$ vanishes on $T_1\otimes T_2\otimes T_3$. This implies that the composition 
$$ V_3\lra  V_1^*\otimes V_2^*\lra  T_1^*\otimes T_2^*$$
has rank at most two, and in fact is equal to two, 
since for $h$ general it drops to rank one in codimension nine,
hence nowhere on the 8-fold $G(2,V_1)\times G(2,V_4)$. The resulting pencil of 
sections of $\cO(1,1)$ on $\PP(T_1)\times \PP(T_2)$ cuts out a degree two subscheme of  
$S_3\subset \PP(V_1)\times \PP(V_2)$, which yields a morphism 
$$\delta : P_H\lra Hilb^2(S_3)=Hilb^2(S).$$

Conversely, a generic point $z$ in $ Hilb^2(S_3)$ is represented by two distinct pairs of lines $(\ell_1,\ell_2)$ and 
 $(\ell_1',\ell_2')$ such that $h$ vanishes on both $\ell_1\otimes \ell_2\otimes V_3$ and $\ell_1' \otimes \ell_2'\otimes V_3$.
Generically the lines $\ell_1$ and $\ell_1'$ are distinct and span a  plane $T_1$, while $\ell_2$ and $\ell_2'$ span a  
plane $T_2$. The morphism $T_1\otimes T_2\lra V_3^*$ induced by $h$ has rank at most two since it
vanishes on $\ell_1\otimes \ell_2$ and $\ell_1'\otimes \ell_2'$, and its image is a two-dimensional
space of linear forms on $V_3$, defining a plane $T_3$. This associates to a generic $z$ in 
$ Hilb^2(S_3)$ a triple $(T_1,T_2,T_3)$  in $G(2,V_1)\times G(2,V_2)\times G(2,V_3)$, which 
belongs to $P_H$. Then it is straightforward to check that this triple is mapped to $z$ by $\delta$. 
Therefore $\delta$ must be a birational morphism. But since  $P_H$ and $Hilb^2(S_3)$ both  have trivial 
canonical bundles, the birationality $\delta$ should be an isomorphism. \qed 

\medskip The fundamental class of $P_H$ in $P = G(2,V_1)\times G(2,V_2)\times G(2,V_3)$,  
as a zero-locus of a section $h_E$ of $E=T_1^*\otimes T_2^*\otimes T_3^*$,  
 is $c_8(T_1^*\otimes T_2^* \otimes T_3^*)$.  
%
The projection of $P_H$ to $G(2,V_1)\times G(2,V_2)$ is an isomorphism with the second degeneracy locus 
of the induced morphism $V_3\lra T_1^*\otimes T_2^*$. By the Thom-Porteous formula, the fundamental class of this  
degeneracy locus is $c_{2,2}(T_1^*\otimes T_2^*)$, 
from where by a standard computation one obtains 
that the 
degree of  $P_H$ with respect to the polarization $\cO(1,1,0)$ is $432$. 

Finally the projection to $G(2,V_1)$ is finite of degree $6$, since a generic fiber can be identified with the 
set of pairs of points inside the intersection of the quartic surface $Q_1\subset \PP(V_1)$ with a generic 
line. Therefore the  polarization  $\cO(1,0,0)$ has degree $12$. 
It has been shown by Ferretti that $P_H$ with
this polarization can be deformed to a double EPW sextic, see \cite{fer}.

\section{The isotropic Grassmannian and its Hilbert scheme}\label{a=1}

In this and the next section we discuss the case $a=1$, 
related to the exceptional group $F_4$. 
The adoint variety $X_{ad}$ for the group $F_4$ has dimension $6a+9 = 15$
and index $3a+5 = 8$, cf. the tables in \S 2 of \cite{cp}.  
Following the procedure from Section \ref{legendrian}, 
we will replace $X_{ad}$ by $X=IG_\omega(3,9)$, the isotropic Grassmannian defined by a 2-form $\omega$ of maximal rank on $V_9$.
This is a Fano variety of the same dimension $15$, and of index $7$. The non semi-simple group $Sp(\omega)\simeq
Sp_9$ acts with two orbits on $X$: 3-planes that contain or do not contain the one-dimensional 
kernel $K$ of $\omega$. 

Following Step 2 in Section \ref{legendrian}, 
the $15$-fold  $X$ is covered by copies of 6-dimensional Legendre cycles $IG(3,6)$ 
(see Table I on p. 89 of \cite{buc}: the case $F_4$, where the Lagrangian Grassmannian 
$IG(3,6)$ of isotropic 3-spaces for a non-degenerate 2-form on $V_6 \cong {\bf C}^6$ is denoted 
by $G_L(3,6)$), and these Legendre cycles are parametrized by an open subset of $P^*=G(6,V_9)$. 
Indeed to any $V_6\subset V_9$ one can associate the variety $X\cap G(3,V_6)\subset G(3,V_9)$.
This is a copy of $IG(3,6)$ when the restriction of $\omega$ to $V_6$ is a non-degenerate 2-form. 

However the restrictions of the general 2-form $\omega$ on $V_9$ to some particular subspaces 
$V_6 \subset V_9$ can be degenerate, which imposes to reconsider the family $P^*$ and to replace 
it with a family $P$ (a birational modification of $P^*$) which will cover the requirements of Step 2 
and the next steps from Section \ref{legendrian}.      
%

\begin{lemm}
The action of  $Sp(\omega)$ on $P^* = G(6,V_9)$ has three orbits $O_i$, $i=0,1,2$, 
defined by the condition that the restriction of $\omega$ to the 6-space $V_6 \in P^*$ 
has rank $6-2i$. The closed orbit $O_2$ has codimension $6$ and is isomorphic 
to $IG(3,8)$.
\end{lemm}

\proof For the general 2-form $\omega$ on $V_9$, the corresponding skew-symmetric map 
$\omega:V_9 \rightarrow V_9^*$ has rank 8. If $K \subset V_9$ is the 1-dimensional kernel of this map, 
one can define the 2-form $\omega$ in terms of a non degenerate 2-form $\omega'$ on an 
eight-dimensional space $V_8$ identified with the image of the map $\omega$, and a 
projection $p: V_9\rightarrow V_8$ from the kernel $K$ of the map $\omega$. Consider some $V_6\subset V_9$, and let 
$2j$  be the rank of the restriction of $\omega$ to $V_6$. If $j=0$, then $p(V_6)\subset V_8$ is isotropic of dimension 
at least five, a contradiction.  If $j=1$, then $p(V_6)$ contains isotropic hyperplanes, which can have dimension at most
four, so $p(V_6)$ has dimension 5, which means that $V_6\supset K$. Moreover by projection, we conclude that 
$O_2$ can be identified with the subvariety of $G(5,V_8)$ defined by the condition that the rank of the restriction 
of $\omega'$ is minimal. This is the minimal orbit of $Sp(8)$ inside $G(5,V_8)\simeq G(3,V_8)$, 
that is  the isotropic Grassmannian $IG(3,8)$. \qed

\medskip\noindent {\it Notation}. We denote by 
$$p: P \rightarrow P^*$$ 
the blowup of $P^* = G(6,V_9)$ along the closed orbit $O_2$, 
by $E = p^{-1}(O_2) \subset P$ the exceptional divisor of $p$, 
and by $\iota : E \rightarrow P$ the embedding of $E$ in $P$.
 
\medskip
The variety $P$ is endowed with a natural vector bundle $F$ of  rank fourteen, 
which we describe below. 

\medskip

Since $p: P - E \rightarrow  P^* - O_2 = O_0 \cup O_1$ is an isomorphism, 
over $P-E$ it is enough to define the fibers of $F$ 
for the 6-spaces $V = V_6 \in O_0 \cup O_1$. 

Denote by $T$ the rank 6 tautological bundle on $P^* = G(6,V_9)$ with fiber $T_V$ over $V = V_6$
identified with the 6-space $V \subset V_9$. For $V \in O_0$ the restriction $\omega_V$ of $\omega$ to $V = T_V$
has rank 6, and by the representation theory of $Sp(6)$ the span of the isotropic grassmannian 
$IG(3,T_V) \cong IG(3,6)$ for $\omega_V$ is the projectivization of the 14-dimensional vector space $F_V$
which fits in the exact sequence
$$0\longrightarrow F_V \longrightarrow \wedge^3T_V \stackrel{\omega_T}{\longrightarrow} T_V \longrightarrow 0,$$ 
where $\omega_T$ denotes the contraction with the 2-form $\omega$. On the orbit $O_1$, although the restriction of 
$\omega$ to $T$ drops rank, the map $\omega_T$ remains surjective and the kernel space $F_V = Ker(\omega_T)$
still has dimension 14.
 
For $V = V_6 \in O_2$,  the restriction $\omega_V$ of $\omega$ to $V = T_V$ has rank two, 
and the image of $\omega_T: \wedge^2 T_V \rightarrow T_V = V$ coincides with the
4-dimensional kernel $T_4\subset V = T_V$ of the  skew-symmetric map $\omega_V: V \rightarrow V^*$.
In particular, for $V \in O_2 \subset P^*$ the rank of the kernel $F_V$ of $\omega_T$ jumps, 
and so the kernel sheaf of $\omega_T$ does not define a vector bundle on $P^*=G(6,V_9)$. 

\medskip


Next, we shall see that after performing a blow-up $P \rightarrow P^*$ 
along the closed orbit $O_2 \subset P^*$, 
the vector bundle $F$ on $P^* - O_2$ can be completed  
to a vector bundle on the projective manifold $P$, 
especially on the exceptional divisor $E$ 
of the blow-up. Due to the natural isomorphism $P - E \cong P^* - O_2$ 
induced by the blow-up, we will keep the notation $F$ also for the 
vector bundle on $P$.   

The normal vector bundle to $O_2$ in $P^*$ is easily seen to 
be isomorphic with $\wedge^2T_4^\vee$. Therefore a point of the exceptional divisor of the blow-up, 
$E\simeq P(N_{O_2/P^*})$, is nothing else than a hyperplane $\Lambda\subset\wedge^2T_4$. 
Now it is straightforward to check that for any such hyperplane one has 
$T_4\wedge\Lambda=\wedge^3T_4$. In particular $T\wedge\Lambda$
is of dimension $14$ inside $\wedge^3T$ 
(if $T_2$ is any complement of $T_4$ in $T$, then $T\wedge\Lambda$ is the direct
sum of  $T_2\wedge\Lambda\simeq  T_2\otimes\Lambda$ of dimension $2\times 5=10$, 
and of  $T_4\wedge\Lambda=\wedge^3T_4$ of dimension $4$). 

\medskip

Define the fiber of $F$ over the point $\Lambda$ of $E$ 
as $F_\Lambda := T\wedge\Lambda$. 
On the exceptional divisor there is also a bundle $p^*T_4$, and we 
denote by $G$ the quotient $p^*(T/T_4)$ which has rank two and is self-dual. 
The above argument yields the following: 

\begin{prop}
$F$ is a rank fourteen sub-bundle of $p^*\wedge^3T $ on $P$.

 There is an exact sequence 
$$0\longrightarrow F \longrightarrow p^*\wedge^3T \stackrel{p^*\omega_T}{\longrightarrow} p^*T 
 \longrightarrow\iota^*G\longrightarrow 0.$$ 
\end{prop}

Let $h$ denote the pull-back to $P$ of the hyperplane class on $P^*$, and $e$ the class of the exceptional divisor $E$. 

\begin{coro} $c_1(F)=-9h+2e$.\end{coro}
 
\begin{coro} The dual bundle $F^\vee$  is globally generated and $H^0(F^\vee) 
\simeq \wedge^3V_9^\vee / \omega\wedge V_9^\vee$. 
\end{coro}

Although we will not really need this, 
it is of worth to mention also the following statement, 
in accordance with our general approach. 

\begin{prop}
The projectivization of any fiber of $F$ cuts $X=IG_\omega(3,9)$ along a generically reduced, 6-dimensional
irreducible subvariety. As a consequence, $P$ parametrizes a family of subvarieties of $X$ made of Legendrian
cycles and degenerations of these. 
\end{prop}

\section{Hyperplane sections}

Next, we follow Steps 3--6 from Section \ref{legendrian} for the case $a = 1$.
 
A hyperplane section of $X_H=X\cap H$ of $X = IG_{\omega}(3,9)$ 
is defined by a 3-form $\Omega$ on $V_9$, 
considered up to $\omega\wedge V_9$.
In this section we study  the family of Legendrian cycles  $IG(3,6)$, and their degenerations, 
that are contained in $H$. 
We denote by  $P_H\subset  P$ the corresponding subvariety. The restrictions of the divisor classes 
$h$ and $e$ from $P$ to $P_H$ will be denoted by $h_H$ and $e_H$.
 
 \begin{prop}
The general 3-form $\Omega$ defines a general section $\tilde\Omega$ of $F^\vee$, 
and the zero-locus of   $\tilde\Omega$ is  $P_H$.
As a consequence, $P_H$ is smooth of dimension $4$ and its canonical class is  
$3e_H$.
\end{prop}

This suggests to study more carefully the intersection of $P_H$ with the exceptional divisor $E$. 
Our main technical result is the following:

\begin{prop}
The intersection $P_H\cap E$ is the disjoint union of nine copies of ${\mathbf P}^3$. 
\end{prop}

\proof Recall that a point $z$ of $E$ is given by a 6-plane $T$ on which $\omega$ has rank two, and a hyperplane 
$\Lambda\subset \wedge^2T_4$, where $T_4\subset T$ is the kernel of the restriction of $\omega$. Let 
$\lambda_4\subset \wedge^2T_4^\vee$ be a linear form defining $\Lambda$. Choose any lift of $\lambda$  to 
 $\wedge^2T^\vee$.  The orthogonal of $F_\Lambda$
in $\wedge^3T^\vee$ is, independently of the latter choice, 
$$F_\Lambda^\perp = \wedge^2T_4^\perp \wedge T^\vee + \lambda\wedge T_4^\perp\subset \wedge^3T^\vee.$$
Next, we choose a basis $e_1, \ldots , e_9$ of $V_9$ adapted to the situation. We will suppose that $T=\langle e_1, \ldots , e_6
\rangle$ and $T_4=\langle e_1, \ldots , e_4\rangle$, so that $T_4^\perp=\langle e_5^*, e_6^*\rangle\subset T^*$
(for simplicity we use the same notation for linear forms on $V_9$ and their restrictions to $T$). 
We can suppose that the restriction 
of $\omega$ to $T$ is $e_5^*\wedge e_6^*$, and then (after changing the basis if necessary), that
$$\omega=e_5^*\wedge e_6^*+e_4^*\wedge e_7^*+e_3^*\wedge e_8^*+e_2^*\wedge e_9^*.$$

Clearly, in this basis the 1-dimensional kernel $K$ of $\omega$ is spanned by $e_1$. 
The condition that $z$ belongs to $P_H$ means that the restriction of $\Omega$ to $T$ belongs to 
$F_\Lambda^\perp$, so that we can write $\Omega$ as 
$$\Omega =e_5^*\wedge e_6^*\wedge u^*+\lambda\wedge v^*+\Omega_7^*\wedge e_7^*
+\Omega_8^*\wedge e_8^*+\Omega_9^*\wedge e_9^*,$$
for some 2-forms $\Omega_7^*, \Omega_8^*, \Omega_9^*$, and some linear forms $u^*$ and $v^*$. Moreover 
$v^*$ is in  $T_4^\perp$, hence it is a linear combination of $ e_5^*$ and $e_6^*$, and we can suppose that $v^*=e_6^*$. 
As a conclusion, 
one is able to write 
the 2-form $\omega$ and the 3-form $\Omega$
as 
$$ \omega=e_5^*\wedge e_6^*+e_4^*\wedge e_7^*+e_3^*\wedge e_8^*+e_2^*\wedge e_9^*,$$
$$ \Omega= \Omega_6^*\wedge  e_6^*+ \Omega_7^*\wedge  e_7^*+\Omega_8^*\wedge  e_8^*+
\Omega_9^*\wedge  e_9^*,$$
for some 2-forms $\Omega_6^*, \Omega_7^*, \Omega_8^*, \Omega_9^*$. 
More intrinsically, 
if the 4-space $W = \langle e_6^*, e_7^*, e_8^*, e_9^* \rangle \subset V_9^*$ 
is the span of $e_6^*,...,e_9^*$, 
then 
$\omega\in V_9^* \wedge W$ 
and  
$\Omega\in \wedge^2 V_9^* \wedge W$. 
Note that if the 5-space $R\in G(5,V_9)$ is the orthogonal to $W$, 
this means exactly that 
$\omega$ and $\Omega$ both vanish identically on $R$.

\begin{lemm} For a generic pair $(\omega,\Omega)$ of a 2-form $\omega$ and a 3-form 
$\Omega$ on the 9-dimensional space $V_9$ there exist
exactly nine 5-dimensional subspaces  $R_1,...,R_9$ of $V_9$ such that  
both $\omega$ and $\Omega$ vanish identically on $R_j$, $j = 1, . . . , 9$. 
\end{lemm}

\proof Denote by $U$ the rank-five tautological vector bundle on the Grassmannian $G(5, V_9)$. The forms 
$\omega$ and $\Omega$ define general sections of the vector bundles $\wedge^2U^*$ and $\wedge^3U^*$ respectively, 
and the common zero locus of these two sections is the set of 5-spaces $R \subset V_9$ 
on which both $\omega$ and $\Omega$ vanish.

Since the two vector bundles $\wedge^2U^*$ and $\wedge^3U^*$ on the 20-dimensional Grassmannian $G(5, V_9)$ are both globally generated and of rank $10$, then 
the zero loci of their generic sections $\omega$ and $\Omega$ 
intersect in a finite number 
$$N=c_{10}(\wedge^2U^*)c_{10}(\wedge^3U^*).$$
of points $R_j \in G(5, V_9), j = 1,...,N$.
It is a standard fact that $c_{10}(\wedge^2U^*)$ is the Schubert class 
$\sigma_{4321}$ defined by the staircase partition. 
Since this Schubert class is self-dual, $N$ is the coefficient of $\sigma_{4321}$ in 
$c_{10}(\wedge^3U^*)$, expressed in terms of Schubert classes. A computation with Sage finally yields 
$N=9$. \qed 

\medskip Let $R_1, \ldots , R_9 \in G(5,V_9)$ be as above, interpreted as 5-spaces in $V_9$. For each $i$, the set of 6-spaces in $V_9$ containing $R_i$ is a 3-dimensional projective space $\pi_i$ linearly 
embedded in $P^* = G(6,V_9)$. Moreover all these $\pi^*_i$ are contained in the orbit 
$O_2$ since (as elements of the zero-locus of $\omega$) the 5-spaces $R_i \subset V_9$ 
are isotropic with respect to $\omega$, see Section \ref{a=1}. 
If $p: P \rightarrow P^* = G(6,V_9)$ is the blowup of $O_2$ as in Section \ref{a=1}, 
and $P_H \subset P$ is the family of Legendrian cycles in $X_H$ as above 
(cf. also with Step 4 in Section \ref{legendrian}),  
then by the preceding:
$$
p(P_H\cap E)\subset \pi_1^*\cup \cdots \cup\pi_9^*.
$$
We claim that each $\pi_i$ can be lifted to a 3-space $\pi_i$ in $E$, and that
$$P_H\cap E=\pi_1\cup \cdots \cup\pi_9.$$
Indeed, once we have written the 3-form $\Omega$ as 
$$
\Omega= \Omega_6^*\wedge  e_6^*+ \Omega_7^*\wedge  e_7^*+\Omega_8^*\wedge  e_8^*+\Omega_9^*\wedge  e_9^*
$$ 
in the adapted basis $e_1,...,e_9$ as in the beginning of this section, 
with corresponding 5-space, say $R_1 = \langle e_1, ..., e_5 \rangle$,  
and its orthogonal 
$W_1 = \langle    e_6^*, e_7^*, e_8^*, e_9^*\rangle$, we get a well-defined map 
$$\bar\Omega : V_9/R_1 \longrightarrow \wedge^2R_1^\vee,$$ 
which is injective for the general pair $(\omega,\Omega)$. 
The map $\bar\Omega$ associates to any 6-space $T$ containing $R_1$ a line in 
$\wedge^2R_1^\vee$, which by restriction is mapped to a line in 
$\wedge^2T_4^\perp$. This defines the lifting
of  $\pi^*_1$ to $\pi_1\subset p^{-1}(\pi_1)\subset E$. \qed

\begin{theo} 
The nine divisors $\pi_1, \ldots , \pi_9$ in $P_H$ can be contracted. The resulting variety 
$Z_H$ is a smooth irreducible  hyperk\"ahler fourfold.
\end{theo}

\proof First observe that the lift $\pi_i$ of $\pi^*_i$ to $E$ being linear, the restriction of $\mathcal{O}_E(1)$
to  $\pi_i$ coincides with the pull-back of $\mathcal{O}_{P^*}(1)$, and also with $\mathcal{O}_{\pi_i}(1)$.
This implies that the canonical class of  $\tilde\pi_i$ is $-4e_{\pi_i}$, and since the canonical class of 
$P_H$ is $3e_H$ we get that the normal bundle of  $\pi_i$ in $P_H$ is 
 $\mathcal{O}_{\pi_i}(-1)$. Therefore the divisor $\pi_i \subset P_H$ can be contracted 
 to a smooth point for any $i = 1,...,9$. Let 
 $$
 f: P_H  \rightarrow Z_H
 $$
 be the composition of the 9 contractions, where by $Z_H$ we have denoted the 
 image variety $Z_H = f(P_H)$. 
 By the preceding, $Z_H$ is a smooth projective variety of dimension $4 = {\rm dim}(P_H)$.  
 By the formulas for blowups, the canonical classes of $P_H$ and $Z_H$ are connected by
the formula
$$3e_H=K_{P_H}=f^*K_{Z_H}+3(\pi_1 +\cdots + \pi_9)=f^*K_{Z_H}+3e_H.$$

Therefore $f^*K_{Z_H}$ is trivial, which implies that the canonical class 
$K_{Z_H}$ by itself is trivial, and as we will see below $Z_H$ is in fact 
a holomorphic symplectic fourfold.  

\medskip

In order to see the last, 
consider a 10-dimensional space $V_{10}$ containing our $V_9$ as a hyperplane, and choose a vector 
$v_0\in V_{10}-V_9$. There is an induced dual decomposition $V_{10}^*=V_9^*\oplus\CC v_0^*$, where 
$v_0^*$ is the linear form evaluating to zero on $V_9$ and to one on $v_0$. Correspondingly, at the level
of three-forms we get the decomposition $\wedge^3V_{10}^*=\wedge^3V_9^*\oplus\wedge^2V_9^*
\wedge v_0^*$. In particular, our 2-form $\omega$ and 3-form $\Omega$ on $V_9$ define a  3-form 
$\Omega_0=\Omega+\omega\wedge v_0^*$ on $V_{10}$. By \cite{dv}, the subvariety $P_{H_0}
\subset G(6,V_{10})$ parametrizing 6-planes on which $\Omega_0$ vanishes identically, is a  
holomorphic symplectic fourfold. The following statement will conclude the proof
of the Theorem.

\begin{prop}
$Z_H$ is isomorphic to  $P_{H_0}$.
\end{prop}

\proof 
The projection of $\PP(V_{10})$ to $\PP(V_9)$ from $\langle v_0\rangle$ induces a rational map 
$$\pi : G(6,V_{10})\lra G(6,V_9).$$ 
This rational map is defined outside the locus $G_0$ of $6$-dimensional subspaces of $V_{10}$ 
containing $v_0$, which is 
a copy of $G(5,V_9)$ and has codimension $4$ in $G(6,V_{10})$. 

Given a $6$-dimensional subspace $T\subset V_9$, a subspace $T_0\subset V_{10}$ such that 
$\pi(T_0)=T$ can be described as the space of vectors $x+u(x)v_0, x\in T$ for some linear form $u$ on 
$T$. Moreover, the fact that $\Omega_0$ vanishes identically on $T_0$ is equivalent to the vanishing of 
$\Omega+u\wedge \omega$ on $T$. This implies that $\pi$ maps $P_{H_0}$ birationally to $P_{H}$. 
We denote the restriction of $\pi$ by $\pi_H : P_{H_0}\lra P_{H}$. 

The rational map $\pi$ can be resolved by just blowing-up $G_0\simeq G(5,V_9)$. Note that $\Omega_0$
vanishes on $T_0=\CC v_0\oplus U$ if and only if $\Omega$ and $\omega$ both vanish on $R$. So 
the (schematic) intersection of $G_0$ with $P_H$ can be identified with the set of nine points
$R_1, \ldots , R_9$  in $G(5,V_9)$ defined in Lemma 2. We can resolve $\pi_H$ by blowing up these
nine points, which yields a birational morphism $\tilde{\pi}_H : \tilde{P}_{H_0}\lra P_{H}$. Moreover the nine 
components of the exceptional divisor are then mapped linearly, hence
isomorphically,  to $\pi_1, \ldots ,\pi_9$. In particular the birational morphism $\tilde\pi_H$ 
preserves the canonical class, so it must be an isomorphism. 
It only remains to contract the nine 
3-spaces on both sides to get an induced isomorphism $\bar{\pi}_H : P_{H_0}\lra Z_{H}$.  
\qed


\bigskip 

\bigskip

{\scriptsize
{\sc Department of Mathematics, Seoul National University, 
Gwanak Campus, Bldg. 27, 
SEOUL 151-747, Korea}

{\it Email address}:  {\tt ailiev@snu.ac.kr}
}

\bigskip

{\scriptsize
{\sc Institut de Math\'ematiques de Marseille,  UMR 7373 CNRS/Aix-Marseille Universit\'e, 
Technop\^ole Ch\^ateau-Gombert, 
39 rue Fr\'ed\'eric Joliot-Curie,
13453 MARSEILLE Cedex 13,
France}

{\it Email address}:  {\tt laurent.manivel@math.cnrs.fr}
}

\end{document}